\documentclass[a4paper,11pt]{amsart}
\usepackage[latin1]{inputenc}
\usepackage[T1]{fontenc}
\usepackage{amssymb,amsfonts,amsmath,amsthm}
\usepackage{color}
\usepackage{setspace}
\usepackage{paralist}
\usepackage{epsfig,graphicx,psfrag}
\usepackage[colorlinks=true, linkcolor=blue, citecolor=red,
urlcolor=blue]{hyperref}

\newcounter{notes}%[page]   %Le 2eme argument fait reinitialiser les numeros de notes a chaque page
\newcommand{\ignore}[1]{}

\newtheorem{theorem}{Theorem}

\newtheorem{conjecture}[theorem]{Conjecture}

\newtheorem{task}[theorem]{Task}

\theoremstyle{definition}
\newtheorem{definition}[theorem]{Definition}

\theoremstyle{remark}
\newtheorem{remark}[theorem]{Remark}

\newtheoremstyle{theoremwithref}{}{}{\itshape}{}{\bfseries}{.}{.5em}{#1 #2 #3}
\theoremstyle{theoremwithref}

\newcommand{\PP}{\mathbf{P}}
\newcommand{\CC}{\mathbf{C}}

\newcommand{\RR}{\mathbf{R}}

\newcommand{\ZZ}{\mathbf{Z}}

\newcommand{\SL}{\mathrm{SL}}
\newcommand{\PSL}{\mathrm{PSL}}
\newcommand{\PGL}{\mathrm{PGL}}

\newcommand{\SO}{\mathrm{SO}}
\newcommand{\PO}{\mathrm{PO}}
\newcommand{\OO}{\mathrm{O}}

\newcommand{\Sp}{\mathrm{Sp}}
\newcommand{\PSp}{\mathrm{PSp}}

\newcommand{\Hom}{\mathrm{Hom}}

\newcommand{\HH}{\mathbb{H}}

\begin{document}

%%%%%%%%%%%%%%%%%%%%%%%%%%%%%%%%%%%%%%%%%%%%%%%%%%%
\title{An invitation to higher Teichm\"uller theory}

%\author[A. Wienhard]{Anna Wienhard}
%\address{Ruprecht-Karls Universit\"at Heidelberg, Mathematisches Institut, Im Neuenheimer Feld~205, 69120 Heidelberg, Germany
%\newline HITS gGmbH, Heidelberg Institute for Theoretical Studies, Schloss-Wolfsbrun\-nenweg 35, 69118 Heidelberg, Germany}
%\email{wienhard@mathi.uni-heidelberg.de}

\thanks{The author acknowledges support by the National Science Foundation under agreement DMS-1536017, by the Sloan Foundation, by the Deutsche Forschungsgemeinschaft, by the European Research Council under ERC-Consolidator grant 614733, and by the Klaus-Tschira-Foundation.}

%\begin{classification}
%\end{classification}

%\begin{keywords}
%Higher Teichm\"uller theory, Total positivity, surface group representations, discrete subgroups of Lie groups, cluster algebras, geodesic flow, mapping class group. 
%\end{keywords}
%
\maketitle

\section{Introduction}
Riemann surfaces are of fundamental importance in many areas of mathematics and theoretical physics. The study of the moduli space of Riemann surfaces of a fixed topological type is intimately related to the study of the Teichm\"uller space of that surface, together with the action of the mapping class group. Classical Teichm\"uller theory has many facets and involves the interplay of various methods from geometry, analysis, dynamics and algebraic geometry. In recent years, higher Teichm\"uller theory emerged as a new field in mathematics. It builds as well on a combination of methods from different areas of mathematics. The goal of this article is to invite the reader to get to know and to get involved into higher Teichm\"uller theory by describing some of its many facets. Along the way we point to open questions, and formulate some conjectures and task for the future. 
We will not be able to discuss every aspect of higher Teichm\"uller theory, and will be very brief on most of them. In particular we will not touch upon universal or infinite higher Teichm\"uller spaces \cite{Labourie_diffeo, Hitchin_diffeo}, or algebraic structures developed in \cite{Labourie_swapping, Sun_swapping, Sun_swapping2}. 
%We will also not have space to discuss the concept of Anosov representations, which provide new classes of discrete subgroups of Lie groups with nice geometric and dynamical properties, but refer the reader to \cite{Kassel_ICM} for this. 

Higher Teichm\"uller theory is concerned with the study of representations of fundamental groups of oriented surface $S$  of negative Euler characteristic into simple real Lie groups $G$ of higher rank. The diversity of the methods involved is due partly to the  one-to-one correspondence between representations, flat bundles, and Higgs bundles given by nonabelian Hodge theory \cite{Simpson_ICM}, which was established in work of Donaldson, Hitchin, Corlette and Simpson, \cite{Donaldson, Donaldson_ASD, Hitchin_selfdual, Corlette88, Simpson}. 

%Given a representation $\rho:\pi_1(S) \to G$ one constructs a flat $G$-bundle on $S$ by taking $E(\rho) = (\widetilde{S} \times G)/\pi_1(S)$, where the action of $\pi_1(S)$ is diagonally, as group of deck transformations on $\widetilde{S}$ and by left multiplication via the representation $\rho$ on $G$. The projection to the first factor give the bundle projection $E(\rho) \to S$, and the product structure $(\widetilde{S} \times G)$ induces a flat connection. Conversely, given a flat $G$ bundle over $S$, the holonomy of the flat connection induces a representation $\pi_1(S) \to G$. 
%When the representation $\rho$ is semisimple, i.e. has reductive Zariski-closure in $G$, there is an $\rho$-equivariant map from $S$ into the symmetric space $X = G/K$. This gives a reduction of the structure group of the flat $G$-bundle to a non-flat $K$-bundle. Choosing a Riemannian metric (or a conformal structure) on $S$, there is a unique harmonic $\rho$-equivariant map from $S$ to $X$. From this one can construct a $G$-Higgs bundle on $S$, which is a holomorphic $G$-bundle together with a holomorphic section of an associated vector bundle. 

We will introduce higher Teichm\"uller spaces below as special subsets of the representation variety $\Hom(\pi_1(S), G))/G$, namely as connected components consisting entirely of discrete and faithful representations. This is however a definition which only arose a posteriori.  
The first family of higher Teichm\"uller spaces, the Hitchin components, has been introduced by Hitchin \cite{Hitchin} using the theory of Higgs bundles. That they are higher Teichm\"uller spaces in the sense of our definition was in general only proven ten years later by Labourie \cite{Labourie_anosov} and independently Fock and Goncharov \cite{FockGoncharov} through the study of the space of positive decorated local systems or positive representations. The second family of higher Teichm\"uller spaces, the spaces of maximal representations, was defined completely independently as the level set of a characteristic number on the representation variety, and its property of being a higher Teichm\"uller space in the sense of our definition was shown by Burger, Iozzi, and Wienhard in \cite{BurgerIozziWienhard_max}, motivated by previous work of Goldman \cite{Goldman_top} in the context of classical Teichm\"uller space. 
The results of \cite{Hitchin}, \cite{Labourie_anosov}, \cite{FockGoncharov} and \cite{BurgerIozziWienhard_max} arose completely independently, from different points of view and using very different methods. 
Only when comparing them it become apparent that the three spaces, Hitchin components, spaces of positive representations, and spaces of maximal representations, have many similarities and provide examples of a new phenomenon. Now we consider them as two families of what we call higher Teichm\"uller spaces. As the reader will see, we are still exploring the similarities and differences of these two families. It is interesting to note that the interplay between geometric and dynamical methods for representations of finitely generated groups and the more analytic and algebro-geometric methods from the theory of Higgs bundles are at the heart of several recent advances in our understanding of higher Teichm\"uller spaces. 

Many questions in higher Teichm\"uller theory are motivated by the things we know about classical Teichm\"uller space, its properties and interesting geometric and dynamical structures it carries. However, there are also several new features that only arise for higher Teichm\"uller spaces and are not present in classical Teichm\"uller theory, see for example Section~\ref{sec:top}, Section~\ref{sec:entropy} and Section~\ref{sec:arithmetic}. 
Higher Teichm\"uller theory is a very young and active field of mathematics. It is shaped by young mathematicians. There are still many open questions and unchartered territory to explore. We therefore hope that many young (and older) mathematicians will accept this invitation and contribute to the field in the future.

\section{Classical Teichm\"uller space} 
Let $S$ be a closed connected oriented topological surface of negative Euler characteristic $\chi(S) = 2-2g < 0$, where $g$ is the genus of $S$. The Teichm\"uller space $\mathcal{T}(S)$ of $S$ is the space of marked conformal classes of Riemannian metrics on $S$. It has been well studied using the theory of quasi-conformal maps as well as methods from hyperbolic geometry. By the uniformization theorem, there is a unique hyperbolic, i.e. constant curvature $-1$, metric in each conformal class. This identifies $\mathcal{T}(S)$ with the moduli space of marked hyperbolic structures. A marked hyperbolic structure is a pair $(X,f_X)$, where $X$ is a hyperbolic surface and $f_X: S \to X$ is an orientation preserving homeomorphism. Two marked hyperbolic structures $(X,f_X)$ and $(Y,f_Y)$ are equivalent if there exists an isometry $g: X \to Y$ such that $g\circ f_X$ is isotopic to $f_Y$. The mapping class group of $S$ acts naturally on $T(S)$ by changing the marking. This action is properly discontinuous, and the quotient of $\mathcal{T}(S)$ by this action is the moduli space $\mathcal{M}(S)$ of Riemann surfaces of topological type given by $S$. Teichm\"uller space is homeomorphic to $\mathbb{R}^{6g-6}$ and the universal cover of $\mathcal{M}(S)$.  

Higher Teichm\"uller theory builds on an algebraic realization of Teichm\"uller space. 
%
%other point of view on Teichm\"uller space through hyperbolic geometry. This point of view has become prominent in the past decades. 
%
%By the uniformization theorem, for any Riemannian metric on $S$, there is a unique hyperbolic, i.e. constant curvature $-1$, metric in its conformal class. Thus we can identify  $\mathcal{T}(S)$ with the moduli space of marked hyperbolic structures. 
%This is formalized by considering pairs  $(X,f_X)$, where $X$ is a hyperbolic surface and $f_X: S \to X$ is an orientation preserving homeomorphism. Two marked hyperbolic structures $(X,f_X)$ and $(Y,f_Y)$ are equivalent if there exists an isometry $g: X \to Y$ such that $g\circ f_X$ is isotopic to $f_Y$. The marking induces a homomorphism $(f_X)_*: \pi_1(S) \to \pi_1(X)$. The mapping class group of $S$ acts naturally on $T(S)$ by changing the marking. This action is properly discontinuous, and the quotient of $\mathcal{T}(S)$ by this action is the moduli space $\mathcal{M}(S)$ of Riemann surfaces of topological type $S$. Teichm\"uller space can hence be identified with the universal cover of $\mathcal{M}(S)$. It is well known that $\mathcal{T}(S)$ is homeomorphic to $\mathbb{R}^{6g-6}$. 
%
%The description of $\mathcal{T}(S)$ through hyperbolic structures on $S$ allows us to give a further, more algebraic realization of Teichm\"uller space, which in fact will be our starting point for the definition of higher Teichm\"uller spaces. 
%
The universal cover $\tilde{X}$ of the hyperbolic surface $X$ naturally identifies with the hyperbolic plane $\HH^2$, and the fundamental group $\pi_1(X)$ acts as group of deck transformations by isometries on  $\tilde{X} \cong \HH^2$. Thus, upon fixing a base point, the marking induces a group homomorphism $(f_X)_*: \pi_1(S) \to \pi_1(X) < \mathrm{Isom}^+(\HH^2) \cong \PSL(2,\RR)$, which is called the holonomy. 
Associating to a marked hyperbolic structure its holonomy gives a well defined injective map 
$$\mathrm{hol}: \mathcal{T}(S) \to \Hom(\pi_1(S), \PSL(2,\RR))/ \PSL(2,\RR).$$ 
The representation variety  $\Hom(\pi_1(S), \PSL(2,\RR))/ \PSL(2,\RR)$ is the space of all group homomorphisms of $\pi_1(S)$ into $\PSL(2,\RR)$, up to conjugation by $\PSL(2,\RR)$. It carries a natural topology (induced from the topology of $\PSL(2,\RR)$).  Teichm\"uller spce $\mathcal{T}(S)$ is a connected component of the representation variety $ \Hom(\pi_1(S), \PSL(2,\RR))/ \PSL(2,\RR)$. It is one of the two connected components, which consist entirely of discrete and faithful representations of $\pi_1(S)$ into $\PSL(2,\RR)$. The other such component is $\mathcal{T}(\overline{S})$, where $\overline{S}$ is the surface with the opposite orientation. 

% 
%We identify Teichm\"uller space with its image under the holonomy map and consider it as a subset  $\mathcal{T}(S) \subset \Hom(\pi_1(S), \PSL(2,\RR))/ \PSL(2,\RR)$. 
%Clearly, any representation in $\mathcal{T}(S)$ is faithful with discrete image. Since $\pi_1(X)$ is a cocompact latttice in $\PSL(2,\RR)$, the set of discrete and faithful representations of $\pi_1(S)$ into $\PSL(2,\RR)$ is a union of connected components. In fact, there are two such components, corresponding to the two possible orientations on $S$, and $\mathcal{T}(S)$ is one of these connected components. 

\begin{remark}
From this point of, as set of discrete and faithful representation, $\mathcal{T}(S)$ was first studied by Fricke. Historically it would thus be more appropriate to call it Fricke space and the generalizations higher Fricke spaces, but it seems hard to change a name that is well established. 
\end{remark}

Classical Teichm\"uller space has many interesting properties and carries additional structure. It is a K\"ahler manifold which admits several Riemannian and non-Riemannian metrics, has nice explicit parametrizations, and carries interesting flows and dynamical systems. We will not be able to recall most of these interesting properties, but will come back to a few of them in the sequel.

%%%%%%Hitchin self-duality ? 

\section{What is higher Teichm\"uller theory?}
We might interpret what higher Teichm\"uller theory is in a narrow or a broader sense. In a very broad sense it is the study of classes of representations of finitely generated groups into Lie groups of higher rank with particularly nice geometric and dynamical behaviour. In the narrow sense one could characterize it as the study of higher Teichmüller spaces as we define them below. In this article we restrict most of our discussion to this narrow interpretation. In the broad sense it is touched upon also in the contributions of Kassel \cite{Kassel_ICM} and Potrie \cite{Potrie_ICM}. 

Teichm\"uller space is a connected component of the representation variety $\Hom(\pi_1(S), \PSL(2,\RR))/ \PSL(2,\RR)$ - this is where higher Teichm\"uller theory takes it starting point. 
Instead of focussing on group homomorphisms of $\pi_1(S)$ into $\PSL(2,\RR)$, we replace $\PSL(2,\RR)$ by a simple Lie group $G$ of higher rank (this is what the higher refers to),  such as $\PSL(n,\RR)$, $n\geq 3$ or $\Sp(2n,\RR)$, $n\geq 2$, and consider the representation variety $\Hom(\pi_1(S), G)/G$. 
We make the following definition: 
\begin{definition}\label{defi:HTS}
A higher Teichm\"uller space is a subset of $\Hom(\pi_1(S), G)/G$, which is a union of connected components that consist entirely of discrete and faithful representations. 
\end{definition}
Note that as soon as $G$ is not locally isomorphic to $\PSL(2,\RR)$,  the group $\pi_1(S)$ is not isomorphic to a lattice in $G$. Therefore the set of discrete and faithful representations is only a closed subset of  $\Hom(\pi_1(S), G)/ G$.  It is thus not clear that higher Teichm\"uller spaces exist at all, and in fact they will only exist for special families of Lie groups $G$. 
In particular, when $G$ is a simply connected complex Lie group, the representation variety $\Hom(\pi_1(S), G)/ G$ is irreducible as an algebraic variety, and hence connected, and there cannot be any connected component consisting entirely of discrete and faithful representations. There are two known families of higher Teichm\"uller spaces, Hitchin components and spaces of maximal representations. They have been discovered from very different points of view and by very different methods. It then became clear that they share many properties, in particular the property requested in Definition~\ref{defi:HTS}. We describe a common underlying characterization, which also suggests the existence of two further families of higher Teichm\"uller spaces in Section~\ref{sec:pos}. 

Hitchin components $\mathcal{T}_H(S,G) $  are defined when $G$ is a split real simple Lie group, 
the space of maximal representations $\mathcal{T}_{max}(S,G)$ is defined when $G$ is a non-compact simple Lie group of Hermitian type. 
In the case when $G = \PSL(2,\RR)$, the Hitchin component and the space of maximal representations agree and coincide with  Teichm\"uller space $\mathcal{T}(S)$. For other groups $G$ not locally isomorphic to $\PSL(2,\RR)$, which are at the same time split and of Hermitian type, i.e. $\Sp(2n,\RR)$ or $\SO(2,3)$, there is a proper inclusion $\mathcal{T}_H(S,G) \subset\mathcal{T}_{max}(S,G)$.

\begin{remark}
We assume that $S$ is a closed surface.There is a related theory for surfaces with punctures or boundary components. However, in this case the corresponding subset of the representation variety  is not a union of connected components. We comment on the situation for surfaces with punctures in Section~\ref{sec:pos}. 
\end{remark}

We shortly review the definitions of Hitchin components and maximal representations. For more details and further properties we refer the reader to the survey \cite{BurgerIozziWienhard_survey}. 
\subsection{Hitchin components}
%Hitchin components were discovered by Hitchin \cite{Hitchin} using methods from the theory of Higgs bundles. This requires the choice of a Riemannian metric on $S$. 
%In the simplest case of $G= \SL(n,\RR)$, a Higgs bundle is a pair $(E,\phi)$ of a holomorphic vector-bundle $E$  over $S$ together with a holomorphic section $\phi \in \Gamma(S,\mathrm{End}(E))$ of the endormorphism bundle. In general there is the theory of $G$-Higgs bundles, and an important correspondence between representations of $\pi_1(S)$ into $G$, flat $G$-bundles on $S$ and $G$-Higgs bundles on $S$. This gives a homeomorphism of moduli spaces 
%$$ 
%\Hom^{ss}(\pi_1(S), G)/G \cong \mathcal{M}_{flat}(S,G) \cong \mathcal{M}_{Higgs}(S,G).
%$$
%Hitchin showed that $\mathcal{M}_{Higgs}(S,G)$ naturally fibers over a space $\mathcal{B}(S,G)$ of $q$-differentials on $S$, where $\mathcal{B}(S,G) = \sum_{i=1}^{k} H^0(S,K^{q_i})$, where $q_i$ are the exponents of the Lie algebra of $G$, and $
%H^0(S,K^{q_i})$ is the space of homolorphic differentials of degree $q_i$ over $S$. 
%In \cite{Hitchin_Lie} Hitchin constructed a Section to the Hitchin fibration $\pi: \mathcal{M}_{Higgs}(S,G) \to \mathcal{B}(S,G) = \sum_{i=1}^{k} H^0(S,K^{q_i})$, the image of which is the Hitchin component. A direct consequence of this is that the Hitchin component is homeomorphic to the vector space $\sum_{i=1}^{k} H^0(S,K^{q_i})\cong \RR^{\mathrm{dim}(G) (2g-2)}$. 
%
%An easier description of the Hitchin components is the following. 
%

Hitchin components are defined when $G$ is a split real simple Lie group. Any split real simple Lie group $G$ contains a three-dimensional principal subgroup, i.e. an embedding $\iota_{p}: \SL(2,\RR) \to G$, which is unique up to conjugation. For the classical Lie groups $\SL(n,\RR)$,  $\Sp(2n,\RR)$, and $\SO(n,n+1)$ this is just the irreducible representation of $\SL(2,\RR)$ in the appropriate dimension. 
Precomposing $\iota_{p}$ with a discrete embedding of $\pi_1(S)$ into $\SL(2,\RR)$ we obtain a representation $\rho_{p}: \pi_1(S) \to G$, which we call a principal Fuchsian representation. 
\begin{definition}\label{defi:Hitchin}
The Hitchin component 
$\mathcal{T}_H(S,G)$ is the connected component of $\Hom(\pi_1(S), G)/G$ containing a principal Fuchsian representation $\rho_{p}:\pi_1(S) \to G$. 
\end{definition}
\begin{remark}
Note that we are a bit sloppy in our terminology, e.g. when $G = \PSL(3,\RR)$ there are $2$ connected components in  $\Hom(\pi_1(S), \PSL(3,\RR))/ \PSL(3,\RR)$ which consists of discrete and faithful representation which preserve the orientation. We refer to each of them as the Hitchin component. 
\end{remark} 

Hitchin showed, using methods from the theory of Higgs bundle that the Hitchin component is homeomorphic to a vector space of dimension $\dim(G) (2g-2)$. He conjectured that these components are geometrically significant and parametrize geometric structures. This was supported by one example. 
Goldman \cite{Goldman_convex} had investigated the spaces of convex real projective structures on $S$ and shown that it is isomorphic to $\RR^{16g-16}$ and that the  holonomy of a convex real projective structure is in the Hitchin component. Soon afterwards Choi and Goldman \cite{ChoiGoldman} proved that in fact $\mathcal{T}_H(S, \PSL(3,\RR))$ parametrizes the space of convex real projective structures on $S$. 
It took another ten years before further progress was made, when Labourie \cite{Labourie_anosov} introduced methods from dynamical systems to the study of representations in the Hitchin component and showed that for $G= \PSL(n,\RR), \PSp(2n,\RR)$, and $\PO(n,n+1)$ representations in the Hitchin component are discrete and faithful. That representations in any Hitchin component are discrete and faithful follows from work of Fock and Goncharov \cite{FockGoncharov}. They investigated the space of positive representations in $\Hom(\pi_1(S), G)/G$, when $G$ is a split real simple Lie group, and showed that it coincides with the Hitchin component (see Section~\ref{sec:pos}).

\subsection{Maximal representations}
Maximal representations are defined when the simple Lie group $G$ is of Hermitian type. They are singled out by a characteristic number, the Toledo number, which for $G = PSL(2,\RR)$ is just the Euler number. 
In \cite{Goldman_top} Goldman showed that the Euler number distinguishes the connected components of $\Hom(\pi_1(S), \PSL(2,\RR))/ \PSL(2,\RR)$, and that Teichm\"uller space corresponds to the connected component formed by representation of  Euler number $2g-2$, which is the maximal value it can attain. In general, the Toledo number is bounded in terms of the Euler characteristic of $S$ and the real rank of $G$, and constant on connected components of $\Hom(\pi_1(S), G)/ G$. 
The space of maximal representations $\mathcal{T}_{max}(S,G)$ is the set of all representations for which the Toledo number assumes it maximal possible value. It is a union of connected components. 
Using methods from bounded cohomology, it was proven in \cite{BurgerIozziWienhard_max} that any maximal representation is faithful with discrete image. 
%
%There only a few simple groups that are at the same time split and of Hermitian type, namely $\PSL(2,\RR), \SO(2,3)$, and $\Sp(2n,\RR)$. 
%For $\PSL(2,\RR)$ Hitchin component and the space of maximal representations agree with $\mathcal{T}(S)$. In  the other cases the Hitchin component $\mathcal{T}_H(S,G)$ is properly contained in the space of maximal representations 
%$\mathcal{T}_{max}(S,G)$. In fact, the space of maximal representations is in general a union of several connected components. We address this in more detail in Section~\ref{sec:top}
%
\begin{remark}
There are two types of Hermitian Lie groups, those of tube type and those not of tube-type. 
Maximal representations into Lie groups that are not of tube type satisfy a rigidity theorem \cite{Toledo_89, Hernandez, BurgerIozziWienhard_max, BradlowGarciaPradaGothen, BradlowGarciaPradaGothen_general}: The image of a maximal representation is always contained in the stabilizer in $G$ of a maximal subsymmetric space of tube type. This reduces the study of maximal representation essentially to the case when $G$ is of tube type. 
\end{remark}

\subsection{Anosov representations} 
Anosov representation are homomorphisms of finitely generated hyperbolic groups $\Gamma$ into arbitrary reductive Lie groups $G$ with special dynamical properties.   
They have been introduced by Labourie \cite{Labourie_anosov} to investigate representations in the Hitchin component, and extended to hyperbolic groups in \cite{GuichardWienhard_anosov}. 
The set of Anosov representations is an open subset of $\Hom(\Gamma, G)/G$, but in general not a union of connected components of $\Hom(\Gamma, G)/G$. Representations in the Hitchin component and maximal representations were the first examples of Anosov representations \cite{Labourie_anosov, BurgerIozziLabourieWienhard, GuichardWienhard_anosov, BurgerIozziWienhard_anosov}. 
We refer to Kassel's contribution \cite{Kassel_ICM} for the definition, more details and more references on Anosov representations. 

The following key properties of Hitchin representations and maximal representations follow from them being Anosov representations (with respect to certain parabolic subgroups). 
\begin{enumerate}
\item Every representation in the Hitchin component and every maximal representation is discrete and faithful. 
\item Let $\rho: \pi_1(S) \to G$ be a Hitchin representation, then there exists a $\rho$-equivariant continuous boundary map $\xi: S^1 \to G/B$ into the generalized flag variety $G/B$, where $B$ is the Borel subgroup of $G$. The map sends distinct points in $S^1$ to transverse points in $G/B$. 
\item Let $\rho: \pi_1(S) \to G$ be a maximal representation, then there exists a $\rho$-equivariant continuous boundary map $\xi: S^1 \to G/S$ into the generalized flag variety $G/S$, where $S$ is a maximal parabolic subgroup of $G$ which fixes a point in the Shilov boundary of the symmetric space $X=G/K$. The map sends distinct points in $S^1$ to transverse points in $G/S$.
\end{enumerate} 

%% and then generalized and further studied by Guichard and Wienhard \cite{GuichardWienhard_anosov}. Further characterizations  Anosov representations provide examples of discrete subgroups in higher rank Lie groups with good dynamical properties. 
%
%
%Let $\rho: \pi_1(S) \to G$ be a Hitchin representation, then 
%\begin{enumerate}
%\item there exists a $\rho$-equivariant continuous boundary map $\xi: S^1 \to G/P$ into the generalized flag variety $G/P$, where $P$ is a minimal parabolic subgroup of $G$. The map sends distinct points in $S^1$ to transverse points in $G/P$.
%\item the representation $\rho$ is faithful and discrete . 
%\item The pull-back of the flat $G$-bundle associated to $\rho$ to the unit tangent bundle $T^1S$ admits a reduction of the structure group to $H$, where $H$ is the Levi subgroup of $P$. 
%\end{enumerate}
%Let $\rho: \pi_1(S) \to G$ be a maximal representation, then 
%\begin{enumerate}
%\item there exists a $\rho$-equivariant continuous boundary map $\xi: S^1 \to G/S$ into the generalized flag variety $G/S$, where $S$ is a maximal parabolic subgroup of $G$ which fixes a point in the Shilov boundary of the symmetric space $G/K$. The map sends distinct points in $S^1$ to transverse points in $G/S$.
%\item the representation $\rho$ is faithful and discrete . 
%\item The pull-back of the flat $G$-bundle associated to $\rho$ to the unit tangent bundle $T^1S$ admits a reduction of the structure group to $L$, where $L$ is the Levi subgroup of $S$. 
%\end{enumerate}
%

\section{Topology of the representation variety}\label{sec:top}
For a connected Lie group $G$ the obstruction to lifting a representation $\pi_1(S) \to G$ to the universal cover of $G$ defines a characteristic invariant in $\mathrm{H}^2(S, \pi_1(G)) \cong \pi_1(G)$. 
For compact simple Lie groups \cite{AtiyahBott} and complex simple Lie groups \cite{Goldman_top, Li_surface}  the connected components of $\Hom(\pi_1(S), G)/G$ 
are in one to one correspondence with elements in $\pi_1(G)$. 
This does not hold anymore for real simple Lie groups in general. Of course, characteristic invariants in $\mathrm{H}^2(S, \pi_1(G))$ still distinguish some of the connected components, but they are not sufficient to distinguish all of them. 
%
%When $G$ is a compact simple Lie group it was shown by Atiyah and Bott \cite{AtiyahBott} that the connected components of $\Hom(\pi_1(S), G)/G$ are in one to one correspondence with elements of $\pi_1(G)$. Consequently, they are all distinguished by the characteristic invariants in $\mathrm{H}^2(S, \pi_1(G))$ which describe the obstruction to lifting a representation $\rho: \pi_1(S) \to G$ to the universal cover of $G$. 
%The analogous statement has been proved by Goldman \cite{Goldmant_top} for $\PSL(2,\RR)$ and $\PSL(2,\CC)$ and conjectured to hold in general for complex Lie groups.  This was then proven by Li \cite{Li_surface}. 
%For real simple Lie groups the situation is however very different. Of course, characteristic invariants in $\mathrm{H}^2(S, \pi_1(G))$ still distinguish some of the connected components, but they are not sufficient to distinguish all of them. 
%We see this in particular when higher Teichm\"uller spaces exist. 
In \cite{Hitchin} Hitchin determined the number of connected components of $\Hom(\pi_1(S), \PSL(n,\RR))/\PSL(n,\RR)$, and showed that Hitchin components have the same characteristic invariants as other components. 
The space of maximal representations, which is defined using characteristic invariants, in fact decomposes itself into several connected components, which hence cannot be distinguished by any characteristic invariant \cite{Gothen, BradlowGarciaPradaGothen_general}. 

A precise count of the number of connected components for several classical groups, and in particular for the connected components of the space of maximal representations has been given using the Morse theoretic methods on the moduli space of Higgs bundles Hitchin introduced \cite{Oliveira, GarciaPradaOliveira_sp, BradlowGarciaPradaGothen_ortho,Gothen, BradlowGarciaPradaGothen_general, GarciaPradaMundet, GarciaPradaGothenMundet}. 
The situation is particularly interesting for $G= \Sp(4,\RR)$ (and similarly for the locally isomorphic group  $\SO^\circ(2,3)$) as there are $2g-4$ connected components in which all representations are Zariski dense \cite{GuichardWienhard_top, BradlowGarciaPradaGothen_SP4}. Maximal representations in these components cannot be obtained by deforming an appropriate Fuchsian representation $\rho: \pi_1(S) \to \SL(2,\RR) \to\Sp(4,\RR)$. An explicit construction of representations in these exceptional components is given in \cite{GuichardWienhard_top}.
Any maximal representations into $\Sp(2n,\RR)$ with $n\geq3$ on the other hand can be deformed either to a principal Fuchsian representation (if it is in a Hitchin component) or to a (twisted) diagonal Fuchsian representation $\rho: \pi_1(S) \to \SL(2,\RR) \times \OO(n) <\Sp(2n,\RR)$. 
In order to distinguish the connected components in the space of maximal representations, additional invariants are necessary. Such additional invariants have been defined on the one hand using methods 
from the theory of Higgs bundles \cite{Gothen, BradlowGarciaPradaGothen_general, Collier_ortho, BaragliaSchaposnik_ortho, A+} and on the other hand using the Anosov property of representations \cite{GuichardWienhard_top}. 
%
%
%
%
%A precise count of the number of connected components of the space of maximal representations has been given for all classical groups  \cite{Gothen, BradlowGarciaPradaGothen_comp}, using the Morse theoretic methods on the moduli space of Higgs bundles Hitchin introduced. 
%The situation is particularly interesting for $G= \Sp(4,\RR)$ (and similarly for the locally isomorphic group  $\SO^\circ(2,3)$) as there are $2g-4$ connected components in which all representations are Zariski dense. This means that these maximal representations can not be obtained by deforming an appropriate Fuchsian representation $\rho: \pi_1(S) \to \SL(2,\RR) \to\Sp(4,\RR)$. An explicit constructions of representations in these coordinates is given in  \cite{GuichardWienhard_top}. 
%For maximal representations into $\Sp(2n,\RR)$ with $n\geq3$ it is in fact the case, that any representation can be deformed either to a principal Fuchsian representation (it it is in a Hitchin component) or to a possible twisted diagonal Fuchsian representation $\rho: \pi_1(S) \to \SL(2,\RR) \times \OO(n) <\Sp(2n,\RR)$. 
%In this case one can also verify that the additional invariants defined by the Anosov property or by a reduction of the Higgs bundle, in fact distinguish all connected components \cite{Gothen, GarciaPradaGothen Mundet, GuichardWienhard_top}.
%We conjecture that this will hold in general. 

We shortly describe the additional invariants arising from the Anosov property. If a representations is Anosov with respect to a parabolic subgroup $P<G$, then the pull-back of the associated flat $G$-bundle to $T^1S$ admits a reduction of the structure group to $L$, where $L$ is the Levi subgroup of $P$. The characteristic invariants of this $L$-bundle provide additional invariants of the representation. 
These additional invariants can be used in particular to further distinguish connected components consisting entirely of Anosov representation. Note that representations can be Anosov with respect to different parabolic subgroups - each such parabolic subgroup gives rise to additional invariants. 
For the case of maximal representations into 
$\Sp(2n,\RR)$ it is shown in \cite{GuichardWienhard_top} that these additional invariants in fact distinguish all connected components. 

\begin{conjecture}\label{conj:top}
Let $G$ be a simple Lie group of higher rank. 
The connected components of $\Hom(\pi_1(S), G)/G$ can be distinguished by characteristic invariants and by additional invariants associated to unions of connected components  consisting entirely of Anosov representations. 
\end{conjecture}

Note that, since Anosov representations are discrete and injective, any  connected component consisting entirely of Anosov representations also provides an example of a higher Teichmüller space. 
Thus Conjecture~\ref{conj:top} implies in particular the following 
\begin{conjecture}\label{conj:top2}
There are connected components of $\Hom(\pi_1(S), G)/G$ which are not distinguished by characteristic invariants, if and only if there exist higher Teichm\"uller spaces in $\Hom(\pi_1(S), G)/G$
\end{conjecture}

Combining Conjecture~\ref{conj:top2} with Conjecture~\ref{conj:pos} gives a precise list of groups (see Theorem~\ref{thm:classification}) for which we expect additional connected components to exist. 
A particularly interesting case is $\SO(p,q)$, $p\neq q$. Here additional components and additional invariants have recently been found via Higgs bundle methods \cite{Collier_ortho, BaragliaSchaposnik_ortho, A+}.

\section{Geometric Structures}\label{sec:geom}
% refer to Fanny's article
Classical Teichm\"uller space $\mathcal{T}(S)$ is not just a space of representations, but in fact a space of geometric structures: every representation is the holonomy of a hyperbolic structure on $S$. For higher Teichm\"uller spaces, such a geometric interpretation is less obvious. The quotient of the symmetric space $Y$ associated to $G$ by $\rho(\pi_1(S))$ is of infinite volume. In order to find geometric structures on compact manifolds associated, other constructions are needed. 

For any representation $\rho: \pi_1(S) \to G$ in the Hitchin component or in the space of maximal representation, there is t a domain of discontinuity in a generalized flag variety $X=G/Q$, on which $\rho$ acts cocompactly. The quotient is a compact manifold $M$ with a locally homogeneous $(G,X)$-structure. This relies on the the construction of domains of discontinuity for Anosov representations given by Guichard and Wienhard in \cite{GuichardWienhard_anosov} and  generalized by Kapovich, Leeb, and Porti in \cite{KLP_dod}. We do not describe this construction here in detail, but refer the reader to \cite{Kassel_ICM}, where locally homogeneous $(G,X)$-structures, Anosov representations and the construction of domains of discontinuity are discussed in more detail. 

The construction of the domains of discontinuity, together with some topological considerations, allows one to deduce the general statement 
\begin{theorem}\cite{GuichardWienhard_anosov} 
For every split real simple Lie group $G$ there exists a generalized flag variety $X$ and a compact manifold $M$ such that $\mathcal{T}_H(S, G)$ parametrizes a connected component of the deformation space of $(G,X)$-structures on $M$. 
For every Lie group of Hermitian type $G$ there exists a generalized flag variety $X$ and a compact manifold $M$ such that for every connected component $C$ of $\mathcal{T}_{max}(S, G)$ the following holds: A Galois cover of $C$ parametrizes a connected component of the deformation space of $(G,X)$-structures on $M$. 
\end{theorem}

In particular, any Hitchin representation or maximal representation is essentially the holonomy of a $(G,X)$-structures on a compact manifold. 
It is however quite hard to get an explicit description of the deformation space of $(G,X)$-structures they parametrize. First it is nontrivial to determine the topology of the quotient manifold, and second it is rather difficult to give a synthetic description of the geometric properties which ensure that the holonomy representation of a $(G,X)$ structure lies in the Hitchin component or in the space of maximal representations. In three cases we such a synthetic description: 
the Hitchin component for $\mathcal{T}_H(S, \PSL(3,\RR))$, the Hitchin component $\mathcal{T}_H(S, \PSL(4,\RR))$ and  $\mathcal{T}_H(S, \PSp(4,\RR))$, the space of maximal representations $\mathcal{T}_{max}(S, \SO^\circ(2,n))$. 
\begin{theorem}\cite{ChoiGoldman}
A representation $\rho: \pi_1(S) \to \PSL(3,\RR)$ is in Hitchin component $\mathcal{T}_H(S, \PSL(3,\RR))$ if and only if it is the holonomy representation of a convex real projective structures on $S$
\end{theorem}
A convex real projective structure on $S$ is a realization of $S$ as the quotient of a convex domain $\Omega \subset \RR\PP^2$ by a group $\Gamma < \PSL(3,\RR)$ of projective linear transformation preserving $\Omega$. 
One aspect which makes this case very special is that the group $\PSL(3,\RR)$ acts as transformation group on the two-dimensional homogeneous space $\RR\PP^2$. The subgroup $\rho(\pi_1(S))$ preserves a convex domain $\Omega$ and acts cocompactly on it. The quotient is a surface homeomorphic to $S$. For more general simple Lie groups of higher rank, there is no two-dimensional generalized flag variety on which they act, and so the quotient manifold $M$ is higher dimensional. 
% refer to convex divisible representations
\begin{theorem}\cite{GuichardWienhard_Duke}\label{thm:convfol}
The Hitchin component $\mathcal{T}_H(S, \PSL(4,\RR))$ parametrizes the space of properly convex foliated projective structures on the unit tangent bundle of $S$. The Hitchin component $\mathcal{T}_H(S, \PSp(4,\RR))$ parametrizes the space of properly convex foliated projective contact structures on the unit tangent bundle of $S$. 
\end{theorem}
\begin{theorem}\cite{CollierTholozanToulisse}
The space of maximal representations $\mathcal{T}_{max}(S, \SO^\circ(2,n))$ parametrizes the space of fibered photon structures on $\OO(n)/\OO(n-2)$ bundles over $S$. 
\end{theorem}
%
%A first step to go beyond these cases it to determine to determine the topology of the compact manifolds that arise as the quotient of a domain of discontinuity. Since the topology is constant on path-connected components of the representation variety, it suffices to determine its topology for one special representation in each connected component. However this turns out to be not an easy task. 
%In each of the above cases $M$ turns out to be a fiber bundle over $S$, and we conjecture this to hold in general: 
\begin{conjecture}[Guichard-Wienhard]\label{conjecture:geom}
Let $\rho: \pi_1(S) \to G$ be a representation in a higher Teichm\"uller space, then there exists a generalized flag variety $X$ and compact fiber bundle $M\to S$, such that $\overline{\rho}: \pi_1(M) \to \pi_1(S) \to G$, where $\pi_1(M) \to \pi_1(S)$ is induced by the bundle map and $\pi_1(S) \to G$ is given by $\rho$, is the holonomy of a locally homogeneous $(G,X)$-structure on $M$. 
\end{conjecture}

In fact, we expect, that for any cocompact domain of discontintuity which is constructed through a balanced thickening in the sense of \cite{KLP_dod}, the quotient manifold is homeomorphic to a compact fiber bundle $M$ over $S$. 
A related conjecture has been made by Dumas and Sanders \cite[Conjecture 1.1]{DumasSanders} for deformations of Hitchin representations in the complexification of $G$. 
They proved the conjecture in the case of $\PSL(3,\CC)$. Guichard and Wienhard \cite{GuichardWienhard_DoDSymp} determine the topology of the quotient manifold for maximal representations and Hitchin representations into the symplectic group,  Alessandrini and Li \cite{AlessandriniLi} prove the conjecture for Hitchin representations into $\PSL(n,\RR)$, and their deformations into $\PSL(n,\CC)$, and Alessandrini, Maloni, and Wienhard \cite{AlessandriniMaloniWienhard} analyze the topology of quotient manifold for complex deformations of symplectic Hitchin representations. 

%
%
%
%
%
%In all the above cases $M$ turns out to be a fiber bundle over $S$, and we conjecture this to hold in general: 
%\begin{conjecture}[Guichard-Wienhard]\label{conjecture:geom}
%Let $\rho: \pi_1(S) \to G$ be a representation in a higher Teichm\"uller space, then there exists a generalized flag variety $X$ and compact fiber bundle $M\to S$, such that $\overline{\rho}: \pi_1(M) \to \pi_1(S) \to G$, where $\pi_1(M) \to \pi_1(S)$ is induced by the bundle map and $\pi_1(S) \to G$ is given by $\rho$, is the holonomy of a locally homogeneous $(G,X)$-structure on $M$. 
%\end{conjecture}
%
%In fact, we even expect, that for any cocompact domain of discontintuity which is constructed by a balanced thickening in the sense of Kapovich-Leeb-Porti, the quotient manifold is homeomorphic to a compact fiber bundle $M$ over $S$. 
%A related conjecture has been made by Dumas-Sanders \cite[Conjecture 1.1]{DumasSanders_1} for deformations of Hitchin representations in the complexification of $G$. 
%They proved the conjecture in the case of $\PSL(3,\RR)$. Guichard-Wienhard\cite{GuichardWienhard_symplmax} show that the conjecture holds for maximal representations into the symplectic group,  Alessandrini-Li \cite{AlessandriniLi} prove it for Hitchin representations into $\PSL(n,\RR), \PSp(2n,\RR)$, and $\SO^\circ(n,n+1)$, and deformations into their complexifications, and Alessandrini-Maloni-Wienhard \cite{AlessandriniMaloniWienhard} analyize the topology of the quotient of the domain of discontinuity for complex deformations of symplectic Hitchin representations. 

It is very interesting to note that the recent advantages \cite{AlessandriniLi, CollierTholozanToulisse} on understanding the topology of the quotient manifolds rely on a finer analysis and description of the Higgs bundle associated to special representations in the Hitchin component or in the space of maximal representations. 
That the explicit description of the Higgs bundles can be used to endow the domain of discontinuity naturally with the structure of a fiber bundle was first described by Baraglia \cite{Baraglia} for $\mathcal{T}_H(S,\PSL(4,\RR))$, where he recovered the projective  structures on the unit tangent bundle from Theorem~\ref{thm:convfol}. 

\section{Relation to the moduli space of Riemann surfaces}
 The outer automorphism group fo $\pi_1(S)$ is isomorphic to the mapping class group of $S$. It acts naturally on $\Hom(\pi_1(S), G)/G$. This action is properly discontinuous on higher Teichm\"uller spaces - in fact more generally on the set of Anosov representations \cite{Labourie_energy, Wienhard_mapping, GuichardWienhard_anosov}. It is natural to ask about the relation between the quotient of higher Teichm\"uller spaces by this action and the moduli space of Riemann surfaces $\mathcal{M}(S)$. For Hitchin components Labourie made a very precise conjecture, based on Hitchin's parametrization of the Hitchin component. We state the parametrization and Labourie's conjecture for $G = \PSL(n,\RR)$ to simplify notation. 
Hitchin introduced the Hitchin component in \cite{Hitchin} using methods from the theory of Higgs bundles. This requires the choice of a conformal structure on $S$. He
Hitchin \cite{Hitchin} 
 showed, using methods from the theory of Higgs bundles, that the Hitchin component is homeomorphic to a vector space. Namely, it is homeomorphic to the space of holomorphic differentials on $S$ with respect to a chosen conformal structure, i.e. $\mathcal{T}_H(S,\PSL(n,\RR))\cong \sum_{i=2}^{n} H^0(S,K^{i})$ This parametrization depends on the choice of a conformal structure and is not invariant under the mapping class group. 
 
 \begin{conjecture}\cite{Labourie_energy}\label{conj:labourie}
 The quotient of $\mathcal{T}_H(S, \PSL(n,\RR))$ by the mapping class group is a holomorphic vector bundle over $\mathcal{M}(S)$, with fiber equal to $\sum_{i=3}^{n} H^0(S,K^{i})$. 
\end{conjecture} 
This conjecture has been proven by Labourie \cite{Labourie_convex} and Loftin \cite{Loftin} for $G = \PSL(3,\RR)$ and by Labourie \cite{Labourie_conj}, using Higgs bundle methods, for all split real Lie groups of rank $2$. It is open in all other cases. 

For maximal representations we expect similarly to get a mapping class group invariant projection from  $\mathcal{T}_{max}(S, G)$ to $\mathcal{T}(S)$, see  \cite[Conjecture~10]{AlessandriniCollier}. 
For $G= \SO(2,n)$ such a projection is constructed in \cite{CollierTholozanToulisse}. For $G$ of rank $2$ Alessandini and Collier  \cite{AlessandriniCollier} construct a mapping class group invariant complex structure on $\mathcal{T}_{max}(S,G)$. For  $G = \Sp(4,\RR)$  they show that the quotient of $\mathcal{T}_{max}(S, G)$ by the mapping class group is a holomorphic vector bundle over $\mathcal{M}(S)$, and describe in detail the fiber over a point, which is rather complicated since the space of maximal representations has nontrivial topology and singular points. 

\section{Positivity}\label{sec:pos}
Hitchin components and maximal representation were introduced and studied by very different methods. It turns out that they do not only share many properties, but also admit a common characterization in terms of positive structures on flag varieties. Only the flag varieties and notions of positivity in question are different for Hitchin components and maximal representations. 

For Hitchin components we consider full flag varieties and Lusztig's total positivity \cite{Lusztig_PosRed}. For maximal representations the flag variety in question is the Shilov boundary of the symmetric space of $G$ and positivity is given by the Maslov cocycle. In order to keep the description simple, we illustrate both notions in examples. We consider $G= \SL(n,\RR)$ for Hitchin components, and $G = \Sp(2n,\RR)$ for maximal representations. 

The relevant flag variety for  $\mathcal{T}_H(S,\SL(n,\RR))$ is the full flag variety
$$
\mathcal{F} (\RR^n):= \{ F=(F_1, F_2, \cdots, F_{n-1}) \, |\, F_i \subset \RR^n, \,  \dim(F_i) = i, \, F_i \subset F_{i+1}\} .
$$
Two flags $F,F'$ are said to be transverse if $ F_i \cap F'_{n-i}  = \{ 0\} $. 
%Given $F$ we denote by $\Omega_F$ the set of all flags in $\mathcal{F}$ which are transverse to $F$. $\Omega_F$ is an open and dense subset of $\mathcal{F}$. 
We fix the standard basis $(e_1, \cdots, e_n)$ of $\RR^n$. 
Let $F \in \mathcal{F}$ be the flag with $F_i = \mathrm{span} (e_1, \cdots, e_i)$, and 
$E \in \mathcal{F}$ the flag with $E_i = \mathrm{span} (e_n, \cdots, e_{n-i+1})$. 

%The group $\GL(n,\RR)$ acts transitively on $\mathcal{F}$ with the subgroups $U$ and $O$ fixing $F$, resp. $E$. 
Any flag $T$ transverse to $F$, is the image of $E$ under a unique unipotent matrix $u_T$. 
The triple of flags $(E, T, F)$ is said to be {\em positive} if and only if $ u_T $ is a totally positive unipotent matrix. Note that a unipotent (here lower triangular) matrix is totally positive if and only of every minor is positive, except those that have to be zero by the condition that the matrix is unipotent. 
Any two transverse flags $(F_1, F_2)$ can be mapped to $(E,F)$ by an element of $\SL(n,\RR)$ and we can extend the notion of positivity to any triple of pairwise transverse flags. 

\begin{theorem}\cite{FockGoncharov, Labourie_anosov, Guichard_CompHit}
Let $\rho: \pi_1(S) \to \SL(n,\RR)$ be a representation. Then $\rho \in \mathcal{T}_H(S, \SL(n,\RR))$ if and only if there exists a continuous $\rho$-equivariant map 
$\xi: S^1 \to \mathcal{F}(\RR^n)$ which sends positive triples in $S^1$ to positive triples in $\mathcal{F}(\RR^n)$.
\end{theorem}

To describe the  analogous characterization of maximal representations into $\Sp(2n,\RR)$ we consider $\RR^{2n}$  with the standard symplectic form $\omega$ and let $\{ e_1, \cdots, e_n, f_1, \cdots f_n\}$ be a symplectic basis. Then  
$$
\mathcal{L} (\RR^{2n}):= \{L \subset \RR^{2n} \, |\, \dim{L} = n, \, \omega|_{L \times L} = 0 \}
$$
is the space of Lagrangian subspaces, and two Lagrangians $L$ and $L'$ are transverse if $L \cap L' = \{ 0\}$. 

Fix 
$L_E = \mathrm{span} (e_1, \cdots, e_n)$ and $L_F = \mathrm{span} (f_1, \cdots, f_n)$. 
Any Lagrangian $L_T \in \mathcal{L}$ transverse to $L_F$ is the image of $L_E$ under an element $v_T =\begin{pmatrix} Id_n& 0\\ M_T& Id_n \end{pmatrix}  \in V$, where $M_T $ is a symmetric matrix.  

The triple of Lagrangians $(L_E, L_T, L_F)$ is said to be {\em  positive} if and only if $M_T \in \mathrm{Pos}(n,\RR) \subset \mathrm{Sym}(n,\RR)$ is positive definite. This is equivalent to the Maslov cocycle of  $(L_E, L_T, L_F)$ being $n$, which is the maximal value it can attain. 
The symplectic group $ \Sp(2n,\RR)$ acts transitively on the space of pairs of transverse Lagrangians and we can extend the notion of positivity to any triple of pairwise transverse Lagrangian. 

 \begin{theorem}\cite{BurgerIozziWienhard_max}
Let $\rho: \pi_1(S) \to \Sp(2n,\RR$ be a representation. Then $\rho \in \mathcal{T}_{max}(S, \Sp(2n,\RR))$ if and only if there exists a continuous $\rho$-equivariant map 
$\xi: S^1 \to \mathcal{L}(\RR^{2n}) $ which sends positive triples in $S^1$ to positive triples in $ \mathcal{L}(\RR^{2n}$. 
\end{theorem}

In \cite{GuichardWienhard_ecm, GuichardWienhard_pos} we introduce the  notion of $\Theta$-positivity. It generalizes Lusztig's total positivity, which is only defined for split real Lie groups, to arbitrary simple Lie groups. There are four families of Lie groups admitting a $\Theta$-positive structure: 
\begin{theorem}\label{thm:classification}\cite[Theorem 4.3.]{GuichardWienhard_ecm}
A simple Lie group $G$ admits a $\Theta$-positive structure if and only if: 
\begin{enumerate}
\item $G$ is a split real form. 
\item $G$ is of Hermitian type of tube type. 
\item $G$ is locally isomorphic to $\SO(p,q)$, $p \neq q$,.
\item $G$ is a real form of $F_4, E_6, E_7, E_8$, whose restricted root system is of type $F_4$. 
\end{enumerate}
\end{theorem}

$\Theta$ refers to a subset of the simple roots. If $G$ admits a $\Theta$-positive structure, then there is  positive semigroup $U^>_\Theta \subset  U_\Theta <  P_\Theta$ with which we can define the notion of positivity for a triple of pairwise transverse points in the generalized flag variety $G/P_\Theta$ as above. Here $P_\Theta$ is the parabolic group associated to the subset of simple roots $\Theta $, and $U_\Theta$ is its unipotent radical.

\begin{definition}\label{defi:thetapos}
A representation $\rho: \pi_1(S)  \to G$ is said to be  $\Theta$-positive if there exists a continuous $\rho$-equivariant map 
$\xi: S^1 \to G/P_\Theta$ which sends positive triples in $S^1$ to positive triples in $ G/P_\Theta$.
\end{definition}

\begin{conjecture}[Guichard-Labourie-Wienhard]\label{conj:pos}
The set of $\Theta$-positive representations $\rho: \pi_1(\Sigma_{g})  \to G$ is open and closed in $ \mathrm{Hom}(\pi_1(S), G)/G$. In particular, $\Theta$-positive representations form higher Teichm\"uller spaces. 
\end{conjecture}
%\begin{conjecture}\label{conjecture:geom}
%Positive representations form higher Teichm\"uller spaces
%\end{conjecture}
For more details on $\Theta$-positivity and $\Theta$-positive representations we refer the reader to  \cite{GuichardWienhard_ecm} and the upcoming papers \cite{GuichardWienhard_pos, GuichardLabourieWienhard}, in which Conjecture~\ref{conjecture:geom} will be partly addressed. In particular we prove that $\Theta$-positive representations are $P_\Theta$-Anosov and form an open subset of $ \mathrm{Hom}(\pi_1(S), G)/G$, and a closed set, at least in the subset of irreducible representations. 
%Recently, Martone and Zhang introduced the notion of positively ratioed representations \cite{MartoneZhang}. We expect $\Theta$-positive representation to be positively ratioed. 
%

The existence of a $\Theta$-positive structure provides a satisfying answer on when and why higher Teichm\"uller spaces exist, and we expect that the families of Lie groups listed in Theorem~\ref{thm:classification} are the only simple Lie groups for which higher Teichm\"uller spaces in $\Hom(\pi_1(S), G)/G$ exist. A particular interesting case is the family of $\Theta$-positive representations for $G= \SO(p,q)$. Here the connected components have recently been determined with Higgs bundle methods, and several of them contain $\Theta$-positive representations  \cite {Collier_ortho, A+}. 

%Key feature: embeddings of SL(2,R) with compact centralizer
%
%

%
%Find common characterization  - Conjecture   
\section{Coordinates and Cluster structures}\label{sec:coordinates}
Teichm\"uller space carries several nice sets of coordinates. The best known are Fenchel-Nielsen coordinates, which encode a hyperbolic structure by the length of and the twist around a set of $3g-3$ disjoint simple closed non-homotopic curves which  give a decomposition of $S$ into a union of $2g-2$ pair of pants. Goldman \cite{Goldman_convex} introduced Fenchel-Nielsen type coordinates on the Hitchin component $\mathcal{T}_H(S,\PSL(3,\RR)$. Here, there are two length and two twist coordinates associated to the curves of a pants decomposition, and in addition two coordinates which associated to each of the pairs of pants. This  a new feature  arises because a convex real projective structure on a pair of pants is not uniquely determined by the holonomies around the boundary. For maximal representations, Fenchel-Nielsen type coordinates were constructed in \cite{Strubel}.

It is often easier to describe coordinates in the situation when the surface is not closed, but has at least one puncture. In this case one can consider decorated flat bundles (or decorated representations), which is a flat bundle or a representation together with additional information around the puncture. Fixing an ideal triangulation, i.e. a triangulation where all the vertices are punctures, this additional information can be used to define coordinates.  Examples of this are Thurston shear coordinates or Penner coordinates for decorated Teichm\"uller space. In the context of higher Teichm\"uller spaces, for decorated representations into split real Lie groups Fock and Goncharov \cite{FockGoncharov} introduced two sets of coordinates, so called $\mathcal{X}$- coordinates, which generalize Thurston shear coordinates,  and $\mathcal{A}$-coordinates, which generalize Penner coordinates. They show that when performing a flip of the triangulation (changing the diagonal in a quadrilateral formed by two adjacent triangles), the change of coordinates is given by a positive rational function.  As a consequence, the set of decorated representations where all coordinates are positive, is independent of the triangulation. In fact, Fock and Goncharov prove that this set of positive representations is precisely the set of positive representations in the sense of Section~\ref{sec:pos}, where the notion of positivity stems from Lusztig's positivity. 
In the case when $G= \PSL(n,\RR)$, the Fock-Goncharov coordinates admit a particularly nice geometric description based on triple ratios and cross ratios. In particular there is a close relation between the coordinates and cluster structures, which received a lot of attention. The change of coordinates associated to a flip of the triangulations is given by a sequence of cluster mutations. This has since been generalized to other classical groups in \cite{Le_1}, see also \cite{Le_2, GoncharovShen} for general split real Lie groups. Related coordinates have been defined by Gaiotto, Moore and Neitzke, using the theory of spectral networks \cite{GaiottoMooreNeitzke, GaiottoMorreNeitzke_snakes}.
For an interpretation of the Weil-Petersson form in terms of cluster algebras see \cite{Gekhtman_etal_weil}.
 
Inspired by Fock-Goncharov coordinates for surfaces with punctures, Bonahon and Dreyer defined coordinates on the Hitchin component $\mathcal{T}_H(S,\PSL(n,\RR))$  and showed that $\mathcal{T}_H(S,\PSL(n,\RR))$  is real analytically homeomorphic to the interior of a convex polygon of dimension $(n^2-1) (2g-2)$ \cite{BonahonDreyer_1, BonahonDreyer_2}. These coordinates are associated to a maximal lamination of the surface $S$ and generalize Thurston's shear coordinates of closed surfaces. A special case for such a maximal lamination is an ideal triangulation of $S$ which is subordinate to a pair of pants decomposition, i.e. the lamination consists of $3g-3$ disjoint simple closed non-homotopic curves which give a pair of pants decomposition, and three curves in each pair of pants, that cut the pair of pants into two ideal triangles. In this case, Zhang \cite{Zhang_Hitchin} provided a reparametrization of the Bonahon-Dreyer coordinates, which give a genuine generalization of Fenchel-Nielsen type coordinates for the Hitchin component  $\mathcal{T}_H(S,\PSL(n,\RR))$, see also \cite{BonahonKim} for a direct comparision with Goldman coordinates when $n=3$. 

In forthcoming work \cite{AlessandriniGuichardRogozinnikovWienhard} we introduce $\mathcal{X}$-type and $\mathcal{A}$-type coordinates for decorated maximal representations of the fundamental group of a punctured surface into the symplectic group $\Sp(2n,\RR)$. These coordinates have the feature of behaving like the coordinates for $G = \PSL(2,\RR)$ but with values in the space of positive definite symmetric matrices $ \mathrm{Pos}(n,\RR)$. In particular, even though they are noncommutative, they exhibit a cluster structure. This structure is similar to the noncommutative cluster structure considered by Berenstein and Retakh \cite{BerensteinRetakh}, except for a difference in some signs. 

It would be interesting to develop similar coordinates for $\Theta$-positive representations, in particular for those into $\SO(p,q)$, and to investigate their properties. The properties of the $\Theta$-positive structure suggests that in this case the cluster-like structure would combine noncommutative and commutative aspects. 
 
\begin{task}\label{prob:cluster}
Develop $\mathcal{X}$-type and $\mathcal{A}$-type coordinates for decorated $\Theta$-positive representations into $\SO(p,q)$. Analyze their cluster-like structures. 
\end{task}

\section{Symplectic geometry and dynamics} 
For any reductive Lie group, the representation variety of a closed surface $\Hom(\pi_1(S), G)/G$ is a symplectic manifold \cite{Goldman_symplectic}. 
On Teichm\"uller space this symplectic structure interacts nicely with Fenchel-Nielsen coordinates. The length and twist coordinates give global Darboux coordinates: the length coordinate associated to a simple closed curve in a pair of pants decomposition is symplectically dual to the twist coordinate associated to this curve, and the symplectic form can be expressed by Wolpert's formula as $\omega = \sum_{i=1}^{3g-3} dl_i \wedge d\tau_i$, where $l_i$ is the length coordinate and $\tau_i$ is the twist coordinate \cite{Wolpert}. 
The twist flows associated to a simple closed curve $c$ on $S$ is the flow given by cutting $S$ along $c$ and  continuously twisting around this curve before gluing the surface back together. 
It is the Hamiltonian flow associated to the length coordinate defined by $c$. The twist flows associated to the $3g-3$ simple closed curves in a pants decomposition on $S$ commute. This gives Teichm\"uller space the structure of a 
complete integrable system. For more general reductive groups, the Hamilitonian flows associated to length functions on $\Hom(\pi_1(S), G)/G$ have been studied by Goldman \cite{Goldman_invariant}. 

In  \cite{SunZhang} provide a new approach to compute the Goldman symplectic form on 
$\mathcal{T}_H(S,\PSL(n,\RR))$. This, in conjunction 
with a companion article by Sun, Zhang and Wienhard \cite{SunWienhardZhang}, gives rise to several nice statements. 
Given maximal lamination with finitely many leaves (and some additional topological data) we construct in \cite{SunWienhardZhang} new families of flows on $\mathcal{T}_H(S,\PSL(n,\RR))$. These flows give a trivialization of the Hitchin component, which is shown to be symplectic in \cite{SunZhang}. Consequently the flows are all Hamiltonian flows and provide $\mathcal{T}_H(S,\PSL(n,\RR)$ the structure of a completely integrable system. 

\begin{task}
The mapping class group $\mathrm{Mod}(S)$ acts naturally on the space of maximal laminations and the additional topological data, so that the symplectic trivializations of $\mathcal{T}_H(S,\PSL(n,\RR))$ induce representations 
$\pi_n: \mathrm{Mod}(S) \to \Sp( \RR^{(n^2-1)(2g-2)})$. Analyze these representations. 
\end{task}
%
% \begin{theorem}\cite{SunWienhardZhang, SunZhang}\label{thm:sympl}
%Fixing an ideal triangulation and additional topological data, there is a symplectic trivialization the Hitchin component $\mathcal{T}_H(S,\PSL(n,\RR))$. 
%\end{theorem}
%This trivialization arises from a construction of new flows on $\mathcal{T}_H(S,\PSL(n,\RR))$ developed in \cite{SunWienhardZhang}, which are in fact all Hamiltonian \cite{SunZhang}. This gives $\mathcal{T}_H(S,\PSL(n,\RR)$ the structure of a completely integrable system. 

A special situation arises when the maximal lamination is an ideal triangulation subordinate to a pants decomposition. 
In this situation we slightly modify the Bonahon-Dreyer coordinates, to get global Darboux coordinates on $\mathcal{T}_H(S,\PSL(n,\RR))$  which consist of $(3g-3) (n-1)$ length coordinates, $(3g-3) (n-1)$ twist coordinates associated to the simple closed curves in the pants decomposition, and $2 \times \frac{(n-1)(n-2)}{2}$ coordinates for each pair of pants. 
The twist flows are the Hamiltonian flows associated to the length functions, and for each pair of pants we introduce  $\frac{(n-1)(n-2)}{2}$ new flows, which we call eruption flows. Their Hamiltonian functions are rather complicated. Nevertheless,  the twist flows and the eruption flows pairwise commute, providing a half dimensional subspace of commuting flows. In the case of $\PSL(3,\RR)$ the eruption flow has been defined in \cite{WienhardZhang}, where it admits a very geometric description.

Classical Teichm\"uller space does not only admit twist flows, but carries several natural flows, for example earthquake flows, which extend twist flows, geodesic flows with respect to the Weil-Petersson metric or the Teichm\"uller metric or, and even an $\SL(2,\RR)$-action. None of this has yet been explored for higher Teichm\"uller spaces. A new approach for lifting Teichm\"uller dynamics to representation varieties for general Lie group $G$ has recently been described by Forni and Goldman \cite{ForniGoldman}.

\section{Geodesic flows and entropy}\label{sec:entropy}
Representations in higher Teichm\"uller space, and more generally Anosov representations, are strongly linked to dynamics on the surface $S$. 
Any such representation gives rise to a H\"older reparametrization of the geodesic flow on $S$ and the representation can essentially be reconstructed from the periods of this reparametrized geodesic flow. 
This dynamical point of view has been first observed by Labourie \cite{Labourie_anosov} and applied by Sambarino \cite{Sambarino_hyperconvex} and has been key in several interesting developments. 

Using the thermodynamical formalism Bridgeman, Canary, Labourie, and Sambarino \cite{BCLS} define the pressure metric on the Hitchin component $\mathcal{T}_H(S,\PSL(n,\RR))$ and more generally spaces of Anosov representations. The pressure metric restricts to a multiple of the Weil-Petersson metric on the subset of principal Fuchsian representations. In the case of $\PSL(3,\RR)$ metrics on the space of convex real projective structures have also been constructed in \cite{GoldmanDarvishzadeh, Li_metric}. 

Other important quantities that have been investigated using this dynamical viewpoint are the critical exponent and the topological entropy of Hitchin representations, which are related to counting orbit points on the symmetric space \cite{Sambarino_counting, Sambarino_quantitative, PollicottSharp}. Here completely new features arise that are not present in classical Teichm\"uller space. On Teichm\"uller space both quantities are constant, but on Hitchin components these functions vary and provide information about the geometry of the representations. There are sequences of representations along which the entropy goes to zero \cite{Zhang_sl3, Zhang_Hitchin}.  The entropy is in fact bounded above. Tholozan for $n=3$ \cite {Tholozan_entropy}, and Potrie and Sambarino  in general, establish an entropy rigidity theorem: A representation saturates the upper bound for the entropy if and only if it is a principal Fuchsian representation \cite{Sambarino_rigidity,Potrie_Sambarino}. This has consequences for the volume of the minimal surface in the symmetric space associated to the representation. A key aspect in the work of Potrie and Sambarino has been the regularity of the map $\xi: S^1 \to G/B$ of a Hitchin representation. 

For maximal representations which are not in the Hitchin component much less is known. One obstacle is the missing regularity of the boundary map $\xi: S^1 \to G/S$, which has rectifiable image, but is in general not smooth. In \cite{GlorieuxMonclair} Glorieux and Monclair study the entropy of Anti-de-Sitter Quasi-Fuchsian representations $\pi_1(S) \to \SO(2,n)$, some of their methods methods might also be useful to investigate maximal representation $\pi_1(S) \to \SO(2,n)$. 

\begin{task}
Investigate the topological entropy of maximal representations. Find and characterize sequences along which the entropy goes to zero. 
Find bounds for the topological entropy and geometrically characterize the representations that saturate these bounds. 
\end{task}

A lot of the geometry of Teichm\"uller space can be recovered from geodesic currents associated to the representations and from their intersection \cite{Bonahon}. Recently Martone and Zhang \cite{MartoneZhang} have associated geodesic currents to positively ratioed representations, a class that includes Hitchin representations and maximal representations but should include also  $\Theta$-positive representations.  In \cite{BCLS_liouville} Bridgeman, Canary, Labourie, and Sambarino define the Liouville current for a Hitchin representation which they use to construct the Liouville pressure metric. 
From the intersection of the geodesic currents one can recover the periods of the reparametrization of the geodesic flow and the periods of crossratios associated to the representation \cite{Labourie_energy, HartnickStrubel}. 
These geodesic currents and the corresponding crossratio functions also play an important role for the next topic we discuss. 

\section{Compactifications}
Classical Teichm\"uller space admits various non-homeomorphic compactifications. One such compactification is the marked length spectrum compactification. 
The marked length spectrum of a representation $\rho:\pi_1(S) \to \PSL(2,\RR)$ associates to any conjugacy class in $[\gamma]$ in $\pi_1(S)$ the translation length of the element $\rho(\gamma)$ in $\HH^2$. 
It is a basic result that for Teichm\"uller space the map $\Phi: \mathcal{T}(S)  \to \PP(\RR^{\mathcal{C}})$ provides an embedding. The closure of $\Phi(\mathcal{T}(S))$ is the marked length spectrum compactification. It is homeomorphic to the Thurston compactification  of $\mathcal{T}(S) $ by the space of projectivized measured laminations. This compactification has been reconstructed using geodesic currents by Bonahon \cite{Bonahon}. 

The marked length spectrum compactification has been generalized to compactifications of spaces of representations of finitely generated groups into reductive Lie groups by Parreau \cite{Parreau}, where the translation length in $\HH^2$ is replaced by the vector valued translation length in the symmetric space associated to $G$. This can be applied to Hitchin components and spaces of maximal representations to provide marked length spectrum compactifications. %(For convex projective structures, the marked length spectrum of has been investigated in \cite{Kim, CooperDelp}.) 
The investigation of the fine structure of these compactifications has just begun. A key ingredient are generalizations of the Collar Lemma from hyperbolic geometry to Hitchin representations \cite{LeeZhang} and to maximal representations \cite{BurgerPozzetti}. See also \cite{LabourieMcShane, VlamisYarmola, FanoniPozzetti} for generalizations of cross-ratio identities in the context of higher Teichm\"uller spaces. 
In \cite{BIPP}  Burger, Iozzi, Parreau and Pozzetti establish a new decomposition theorem for geodesic currents, which will play a crucial role in their program to understand the marked length spectrum compactification of maximal representations. 

Compactifications of the space of positive representations (when $S$ has punctures) have been constructed using explicit parametrizations and the theory of  tropicalizations \cite{FockGoncharov, Alessandrini, Le_lamination}. 

In all these constructions it is a challenge to give a geometric interpretation of points in the boundary of the compactification. The most natural is in terms of actions on $\RR$-buildings \cite{Parreau, Le_lamination, BurgerPozzetti}. This naturally generalizes the description of boundary points of Teichm\"uller space by actions on $\RR$-trees. However, Thurston's compactification gives an interpretation of boundary points by measured laminations on $S$. It would be interesting to get a description of boundary points of Hitchin components or spaces of maximal representations in terms of geometric objects on $S$ that generalize measured laminations, and to relate the compactifications to degenerations of the geometric structures associated to Hitchin components and maximal representations in low dimensions (see Section~\ref{sec:geom}). For convex real projective structures such an interpretation of boundary points in terms of a mixed structure consisting of a measured lamination and a HEX-metric has been announced by Cooper, Delp, Long, and Thistlethwaite \cite{CDLT}.

\section{Arithmetics}\label{sec:arithmetic}.
Discrete Zariski dense subgroups of semisimple Lie groups have been studied for a long time. In recent years there has been a revived interest from number theory in discrete Zariski dense subgroups, which are contained in arithmetic lattices without being lattices themselves. Such groups have been coined thin groups \cite{Sarnak}. Hitchin representations and maximal representations, and more generally Anosov representations, provide many examples of discrete Zariski dense subgroups of higher rank Lie groups which are not lattices. They are not necessarily contained in (arithmetic) lattices. However, there are many Hitchin representations and maximal representations that are even integral, i.e. up to conjugation contained in the integral points of the group $G$. This is a new feature that arises for higher Teichm\"uller spaces and is not present in classical Teichm\"uller space. 
The group $\PSL(2,\ZZ)$ is a non-uniform lattice in $\PSL(2,\RR)$, consequently there is no discrete embedding of the fundamental group $\pi_1(S)$ of a closed oriented surface into $\PSL(2,\RR)$ which takes values in the integer points $\PSL(2,\ZZ)$.%There are however plenty embeddings of  $\pi_1(S)$ into the integer points $\PSL(n,\ZZ)$ which lie in the Hitchin component $\mathcal{T}_H(\pi_1(S), \PSL(n,\RR))$ for $n$ odd and $n\geq 3$. 

For $n=3$ the first examples of integral Hitchin representations were given by Kac and Vinberg \cite{KacVinberg}, infinite families were constructed by Long, Reid, and Thistlethwaite with an explicit description using triangle groups. 
\begin{theorem}\label{thm:LRT}\cite{LongReidThistlethwaite}
Let $\Gamma = \langle a, b \, |\, a^3 = b^3 = (ab)^4 = 1 \rangle$ be the (3,3,4) triangle group. 
There is an  explicit polynomial map $\phi: \RR \to \Hom(\Gamma, \PSL(3,\RR))$ whose image lies in the Hitchin component. For all $t\in \ZZ$,  the image of $\phi(t)$ give a Zariski-dense subgroup of $\PSL(3,\ZZ)$. The representations $\phi(t)$, $t\in Z$ are pairwise not conjugate in $\PGL(3,\RR)$. 
\end{theorem}
Since $\Gamma$ contains subgroups of finite index which are isomorphic to the fundamental group $\pi_1(S)$ of a closed oriented surface $S$, Theorem~\ref{thm:LRT} gives rise to infinitely many, non conjugate integral representations in $\mathcal{T}_H(\pi_1(S), \PSL(n,\RR))$. The representations lie on different mapping class group orbits. 

In unpublished work with Burger and Labourie we use bending to show the following 
\begin{theorem}\label{thm:BLW}
For $n\geq 5$ and odd there are infinitely many pairwise non conjugate integral representations in the Hitchin component $ \mathcal{T}_H(\pi_1(S), \PSL(n,\RR))$. These representations lie on different mapping class group orbits. 
\end{theorem}

\begin{task}
Develop tools to count integral representations in $\mathcal{T}_H(\pi_1(S), \PSL(n,\RR))$ modulo the action of the mapping class group. Investigate the counting functions and their asymptotics. 
\end{task}

A first step to start counting is to find appropriate height functions on $\mathcal{T}_H(\pi_1(S), \PSL(n,\RR))$ such that there are only finitely many integral representations of finite height. A height function, inspired by Thurston's asymmetric metric on Teichm\"uller space has been proposed by Burger and Labourie. 

%%%%%% Chekc Marc's talk - Tengren's result. 

For more general constructions of surface subgroups in lattices of Lie groups, following the construction of Kahn and Markovic surface subgroups in three-manifold groups \cite{KahnMarkovic} we refer to work of Hamenstädt and Kahn \cite{HamenstaedtKahn}, and forthcoming work of Labourie, Kahn, and Mozes \cite{LabourieKahnMozes}. Examples of  integral maximal representations are constructed in \cite{Toledo_hilbert}.

\section{Complex Analytic Theory}
In classical Teichm\"uller theory complex analytic methods and the theory of quasi-conformal mappings play a crucial role. These aspects are so far largely absent from higher Teichm\"uller theory. 
Dumas and Sanders started exploring the complex analytic aspects of discrete subgroups of complex Lie groups of higher rank in \cite{DumasSanders}. They investigate in particular deformations of Hitchin representations and maximal representations in the complexifications of $G$, and establish important properties of the complex compact manifolds $M$ that arise as quotients of domains of discontinuity of these representations (see Section~\ref{sec:geom}). In a forthcoming paper \cite{DumasSanders_2} the complex deformation theory of these representations will be investigated further. 

\section{Higher dimensional higher Teichm\"uller spaces} 
Fundamental groups of surfaces are not the only finitely generated groups for which there are special connected components in the representation variety $\Hom(\pi_1(S), G)/G$, which consist entirely of discrete and faithful representations. This phenomenon also arises for fundamental groups of higher dimensional manifolds, and even for more general finitely generated hyperbolic groups. 
The main examples are convex divisible representations, which have been introduced and studied by Benoist in a series of papers, starting with \cite{Benoist_CD1}. They are generalizations of convex real projective structures on surfaces, and exist in any dimension. 

 Let $N$ be a compact hyperbolic manifold of dimension $n$ and $\pi_1(N)$ its fundamental group. 
A representation $\rho:\pi_1(N) \to \PGL(n+1,\RR)$ is convex divisible if there exists a strictly $\rho(\pi_1(N))$- invariant convex domain in $\RR\PP^n$, on which $\rho(\pi_1(N))$ acts cocompactly. 
\begin{theorem}\cite{Benoist_convex_III}
The set of convex divisible representations is a union of connected components of $\Hom(\pi_1(N), \PGL(n+1,\RR))/\PGL(n+1,\RR)$ consisting entirely of discrete and faithful representations. 
\end{theorem}

Barbot \cite{Barbot} showed that a similar phenomenon arise for representations $\rho:\pi_1(N) \to \SO(2,n)$ that are Anti-de-Sitter Quasi-Fuchsian representation. In fact these are precisely Anosov representations of $\pi_1(N)$ into $\SO(2,n)$ with respect to the parabolic subgroup that stabilizes an isotropic line \cite{BarbotMerigot}. 

\begin{theorem}\cite{Barbot}
The set of Quasi Fuchsian AdS representations in $\Hom(\pi_1(N), \SO(2,n))/\SO(2,n)$ is a connected component consisting entirely of discrete and faithful representations.
\end{theorem}

In view of Definition~\ref{defi:HTS} we might call these connected components of the representation variety $\Hom(\pi_1(N), G)/G$ containing only discrete and faithful representations higher dimensional higher Teichm\"uller spaces. 

When $N$ is of dimension two, the notion of $\Theta$-positivity gives us a conjectural criterion why and when higher Teichm\"uller spaces exist. It would be interesting to find a unifying principle behind the existence of such special connected components in $\Hom(\pi_1(N), G)/G$.  
%A first attempt could be to analyze the connected components containing standard representations $\rho_0: \pi_1(N) \to \SO(1,n) \to G$, where $\pi_1(N) \to \SO(1,n)$ is the unique realization of $\pi_1(N)$ as a lattice. 

\begin{task}
Find the underlying principle for the existence of connected components of the representation variety $\Hom(\pi_1(N), G)/G$ which consist entirely of discrete and faithful representations. 
\end{task}

A first test case could be to analyze deformations of the representation $\rho_0: \pi_1(N) \to \SO(1,n) \to \SO(k,n)$,  or more generally deformations of representations $\rho_0: \pi_1(N) \to \SO(1,n) \to G$, where the centralizer of $\SO(1,n)$ in $G$ is compact and  contained in the maximal compact subgroup of the Levi group of a parabolic subgroup containing the parabolic subgroup defined by $\SO(1,n)$.  We expect all such deformations to be discrete and faithful.

\end{document}